\documentclass{amsart}

\usepackage{amsthm}
\usepackage{graphics}
\usepackage{epsfig}
\theoremstyle{plain}
\newtheorem{theorem}{Theorem}[section]
\newtheorem{proposition}{Proposition}[section]

\theoremstyle{definition}
\newtheorem{definition}{Definition}[section]
\theoremstyle{remark}
\newtheorem{remark}{Remark}[section]
\newtheorem{example}{Example}[section]

\newcommand{\refE}[1]{(\ref{E:#1})}
\newcommand{\refS}[1]{Section~\ref{S:#1}}

\renewcommand{\a}{\ensuremath{\alpha}}

\newcommand{\C}{\ensuremath{\mathbb{C}}}

\renewcommand{\P}{\ensuremath{\mathbb{P}}}
\newcommand{\Z}{\ensuremath{\mathbb{Z}}}

\newcommand{\cins}{\frac 1{2\pi\mathrm{i}}\int_{C_S}}
\newcommand{\A}{\mathcal{A}}
\newcommand{\tr}{\mathrm{tr}}

\newcommand{\g}{\mathfrak{g}}
\newcommand{\gb}{\overline{\mathfrak{g}}}
\newcommand{\gh}{\widehat{\mathfrak{g}}}

\newcommand{\Vh}{\widehat V}

\renewcommand{\L}{\mathcal{L}}
\newcommand{\ord}{\operatorname{ord}}
\newcommand{\res}{\operatorname{res}}
\newcommand{\ldot}{\,.\,}
\newcommand{\de}{\delta}

\newcommand{\Ho}{\mathrm{H}}
\newcommand{\lpz}{\C[z,z^{-1}]}
\newcommand{\sln}{\mathfrak{sl}}

\newcommand{\gl}{\mathfrak{gl}}
\begin{document}

\vspace*{-1cm}
\hbox{ }
{{\hspace*{\fill} math.QA/0510440}}

\vspace*{2cm}

\title{Higher Genus Affine Lie Algebras 
\\
of Krichever -- Novikov Type}

\author{MARTIN SCHLICHENMAIER}

\address{
University of Luxembourg, \\
Campus Limpertsberg, \\
162A, Avenue de la Faiencerie, \\
L-1511 Luxembourg, Grand-Duchy of Luxembourg \\
E-mail: martin.schlichenmaier@uni.lu}

\begin{abstract}
Classical affine Lie algebras appear e.g. as symmetries of
infinite dimensional integrable systems and
are related to certain differential equations.
They  are central extensions of current
algebras associated to finite-dimensional Lie algebras $\mathfrak{g}$.
In geometric terms these current algebras might be described as
Lie algebra valued meromorphic functions on the
Riemann sphere with two possible poles. They carry a natural grading.
In this talk the generalization to higher genus 
compact Riemann surfaces and more
poles is reviewed.
In  case that the Lie algebra $\mathfrak{g}$ is reductive
(e.g. $\mathfrak{g}$ is simple, semi-simple, abelian, ...)
a complete classification of (almost-) graded central extensions is
given.
In particular, for $\mathfrak{g}$ simple  there exists a unique
non-trivial (almost-)graded extension class.
The considered algebras are related to difference 
equations, special functions and play  a role in
Conformal Field Theory.
\end{abstract}
\date{19.10.2005}
\maketitle

\begin{minipage}[c]{0.9\linewidth}
{\it Talk presented at the International Conference on Difference
Equations, Special Functions, and Applications, Munich, July
2005}
\end{minipage}

\vskip 1cm
\section{Introduction}
%%%%%%%%%%%%%%%%%%%%%%%%%%%%%%%%%%%%%%%%%%%%
Classical current algebras (also called loop algebras) and their 
central extensions, the affine Lie algebras, are of fundamental 
importance in quite a number of fields in mathematics and its
applications. These algebras are examples of infinite dimensional 
Lie algebras which are still tractable. They constitute the subclass of
Kac-Moody algebras of untwisted affine type.

If one rewrites the original purely algebraic definition in geometric terms the
classical current algebras correspond to Lie algebra valued 
meromorphic functions on the Riemann sphere (i.e.~on the
unique compact Riemann surface of genus zero) which are allowed to
have poles only at two fixed points.

If this rewriting is done, a very useful generalization 
(e.g.~needed in string theory) is to consider the 
objects over a compact Riemann surface of arbitrary genus with more than two
points where poles are allowed.
The main problem is to introduce a replacement of the grading in the classical
case, which is necessary to construct highest weight and Fock space
representations.
This is obtained by an almost-graded structure (see \refS{hg}), a weaker 
structure but still strong enough to do the job.
Furthermore, to obtain representations of certain types one is forced to 
pass over to central extensions.

Such objects (vector fields, functions, etc.) and central extensions
for higher genus with two possible poles were introduced by Krichever and Novikov
\cite{KNFa} and generalized by me to the 
multi-point situation \cite{Schlce}.
These objects are of importance in a global operator approach to
Conformal Field Theory \cite{SchlShwznw},\cite{SchlSh}.
More generally, the current algebra resp. 
their central extensions, the affine algebras, correspond to
symmetries of infinite dimensional systems.
Their $q$-deformed version (i.e.~the quantum affine algebras) are in close
connection with difference equation and with special functions.

There is  a very interesting direct relation to difference equation.
Krichever and Novikov constructed higher genus analogues of Baker -- Akhiezer 
difference functions which are eigen functions of suitable difference
equations.
Starting from 
these functions, representations of higher genus affine algebras associated
to the Lie algebra  $\sln(2,\C)$ are obtained. These are related to the 
two-dimensional Toda lattice \cite{Ktoda}.

In this write-up of the talk I report on uniqueness and classification results
for
higher genus multi-point affine Lie algebras which I recently obtained.
In particular, it turns out that for the current algebra associated to a 
finite-dimensional simple Lie algebra (e.g. for $\sln(n,\C)$) there exists 
up to equivalence and rescaling a unique
non-trivial central extension which extends the almost-grading of the
current algebra.
The proofs can be found in \cite{Schlaff}. There also further references,
historical remarks and corresponding results for the Lie algebras of Lie
algebra
valued meromorphic differential  operators can be found.
The results depend also on a complete classification of central
extensions of scalar functions, vector fields and differential operators
of Krichever--Novikov type
obtained in \cite{Schlloc}.

I am indebted to Paul Terwilliger 
who asked me to supply explicit examples of 
such algebras. They can be found in \refS{zero} and \refS{torus}.

%%%%%%%%%%%%%%%%%%%%%%%%%%%%%%%%%%%%%%%%%%%%%%%
\section{The classical situation and some algebraic background}
%%%%%%%%%%%%%%%%%%%%%%%%%%%%%%%%%%%%%%%%%%%%%%%%
Let us first consider the nowadays classical affine Lie algebras.
Let $\g$ be a finite-dimensional complex Lie algebra.
A special example of fundamental importance is given by the
algebra of trace-less matrices
\begin{equation}
\sln(n,\C):=\{A\in Mat(n,\C)\mid \tr(A)=0\},
\end{equation}
with $[A,B]:=AB-BA$ the commutator as Lie product.

The current algebra $\gb$
(sometimes also called loop algebra) is obtained
by tensoring $\g$ by the (associative and commutative)
algebra $\lpz$ of Laurent polynomials, i.e. 
$\gb=\g\otimes \lpz$ with the Lie product
\begin{equation}\label{E:classalg}
[x\otimes z^n,y\otimes z^m]:=[x,y]\otimes z^{n+m},
\quad x,y\in \g, \ n,m\in\Z.
\end{equation}
If $\g$ is a matrix algebra, then $\gb$ can be considered as matrices
with Laurent polynomials as entries, e.g. 
a typical element of $\overline{\sln}(2,\C)$ can be written as
\begin{equation}
\begin{pmatrix}
a(z,z^{-1})&\quad& b(z,z^{-1})
\\
c(z,z^{-1})&\quad &  - a(z,z^{-1})
\end{pmatrix},
\end{equation}
where $a,b,c$ are polynomials in $z$ and $z^{-1}$.

By setting $\deg(x\otimes z^n):=n$ the Lie algebra $\gb$ is
graded (see \refE{classalg}).
Clearly, $\gb$ is an infinite dimensional Lie algebra.
These algebras  
are candidates for symmetry algebras of systems with infinitely many 
independent symmetries.
Unfortunately, in the process of constructing representations of
certain types (e.g. highest weight representations) one is forced
to ``regularize'' certain natural actions. As a result one 
obtains only so called {\it projective representations}, 
which in turn define 
honest representations of certain central extensions of $\gb$.

What is a central extension $\Vh$ of a Lie algebra $V$?
As vector space we take $\Vh=V\oplus \C$.
We set $t:=(0,1)$ and $\hat a:=(a,0)$ and consider the following
product
\begin{equation}\label{E:central}
[\hat a,\hat b]:=
\widehat{[a,b]}+\psi(a,b)\,t,\quad
[\hat a,t]=0,
\end{equation}
with $\psi:V\times V\to V$ a bilinear map.
Now $\Vh$ is a Lie algebra if and only if $\psi$ is a Lie algebra
two-cocycle of $V$ (with values in the trivial module $\C$).
The cocycle conditions are 
\begin{equation}
\psi(a,b)=-\psi(b,a),
\qquad
\psi([a,b],c)+
\psi([b,c],a)+
\psi([c,a],b)=0,
\end{equation}
for all $a,b,c\in V$.

Two central extensions of the same algebra are called 
equivalent if they are the same up to some change of basis
of the type $\hat a\mapsto \tilde a=(a,\phi(a))$.
In more precise terms, given two extensions defined by
$\psi_1$ and $\psi_2$ respectively,  the two extensions
are called equivalent if there exists a linear form 
$\phi:V\to \C$ such that
\begin{equation}
\psi_1(a,b)-
\psi_2(a,b)=\phi([a,b]).
\end{equation}
In other words, the difference is a Lie algebra cohomology coboundary.

\medskip

{\bf Fact.}
The set of equivalence classes of central extensions is via
\refE{central}
in
1 to 1 correspondence to the space of Lie algebra two-cohomology classes
$\Ho^2(V,\C)=Z^2(V,\C)/B^2(V,\C)$
(cocycles modulo coboundaries).

\medskip
How do we obtain central extensions for our current algebra?
Let  $\a$ be an invariant, symmetric bilinear form for $\g$.
Invariance means that 
$\alpha([a,b],c)=\alpha(a,[b,c])$ for all $a,b,c\in\g$.
For a simple Lie algebra the Cartan-Killing form is up to a
rescaling the only such form. In particular for $\sln(n,\C)$ 
it is given by 
$ \alpha(A,B)=\tr(AB)$.
Then  a central extension $\gh=\gb\oplus\C t$
is defined by
\begin{equation}\label{E:centralclass}
[x\otimes z^n,y\otimes z^m]=[x,y]\otimes z^{n+m}-\alpha(x,y)
\cdot n\cdot \de_{n}^{-m}\cdot t.
\end{equation}
To avoid cumbersome notation I dropped the $\ \hat{ }\  $ in the notation.
It is called the (classical) affine Lie algebra associated to $\g$.
By setting $\deg t:=0$ (and using $n=\deg(x\otimes z^n)$) the 
affine algebra is a graded algebra.
If the finite dimensional Lie algebra $\g$ is simple then 
$\gh$ defined via \refE{centralclass} is up to equivalence 
of extensions and
rescaling of the central element the only non-trivial
central extension of $\gb$.
In this case the algebras are exactly the Kac-Moody algebras of
untwisted affine type \cite{Kac}.

%%%%%%%%%%%%%%%%%%%%%%%%%%%%%%%%%%%%%%%%%%%%%%%
\section{The  higher genus case}\label{S:hg}
%%%%%%%%%%%%%%%%%%%%%%%%%%%%%%%%%%%%%%%%%%%%%%%
Before we can extend the construction to higher genus we have to
geometrize the classical situation.
Recall that the associative algebra of Laurent polynomials $\lpz$ 
can equivalently be described as the algebra consisting of those 
meromorphic functions on the Riemann sphere (resp. the complex
projective line $\P^1(\C)$) which are holomorphic outside
$z=0$ and $z=\infty$ ($z$ the quasi-global coordinate).
The current algebra $\gb$ can be interpreted as Lie algebra of
$\g$--valued
meromorphic functions on the Riemann sphere with possible poles
only at $z=0$ and $z=\infty$.

The Riemann sphere is the unique compact Riemann surface of
genus zero. From this point of view the next step is to take
$X$ any compact Riemann surface of arbitrary genus $g$ and
 an arbitrary finite set $A$ of points where poles 
of the meromorphic objects will be allowed.
In this way we obtain the higher genus (multi-point) current algebra
as the algebra of $\g$--valued functions on $X$ with only possibly 
poles at $A$.
We need gradings, central extensions etc.

\begin{figure}
\begin{minipage}[b]{\linewidth}
\centering\epsfig{figure=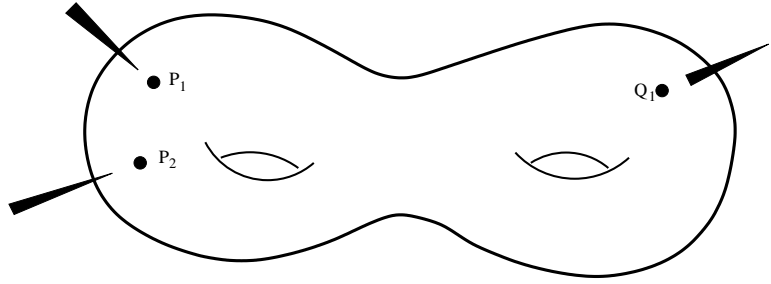,width=.8\linewidth,clip=}
\caption{The higher genus, multi-point case}
\label{figg2}
\end{minipage}\hfill
\end{figure}

For this goal we split $A$ into two non-empty disjoint subsets
$I$ and $O$, $A=I\cup O$.
In the interpretation of string theory, $I$ corresponds to incoming
free strings and $O$ to outgoing free strings.
Let $K$ be the number of points in $I$.
See Figure \ref{figg2} for an example given by a 
 Riemann surface of genus two with
$I=\{P_1,P_2\}$ and $O=\{Q_1\}$.
%%%%%%%%%%%%%%%%%%%%%%%%%%%%%%%%%%%%%%%%%%%%%%%

%%%%%%%%%%%%%%%%%%%%%%%%%%%%%%%%%%%%%%%%%%%% 
Let $\A$ be the associative algebra of functions 
meromorphic on $X$ and holomorphic outside of $A$.
In some earlier work \cite{Schlce} I introduced 
\begin{equation}
\{A_{n,p}\mid  n\in\Z,\ p=1,\ldots,K\}
\end{equation}
a certain adapted basis of $\A$.
For the exact definition I refer to this publication.
Here we only note that
\begin{equation}
\ord_{P_i}(A_{n,p})=n+1-\de_{i}^p,\quad \forall P_i\in I.
\end{equation}
For genus zero and $I=\{0\}$, $O=\{\infty\}$
we get $A_{n,p}=z^n$.
Let $\A_n:=\langle A_{n,p}\mid p=1,\ldots, K\rangle$
be the $K$-dimensional subspace of $\A$. We have
$\A=\oplus_{n\in\Z} \A_n$ and there exist constants
$L_1,L_2$ (independent of $n$ and $m$) such that
\begin{equation}
\A_n\cdot A_m\subseteq \bigoplus_{h=n+m-L_1}^{n+m+L_2} \A_h,
\qquad  \forall n,m\in\Z.
\end{equation}
We call the elements of $\A_n$ homogeneous elements of degree $n$.
As long as $L_1$ and $L_2$ cannot to be chosen to be 0 
the algebra is  not honestly graded. It is only {\it almost-graded}.
In a similar way one introduces almost-gradedness for Lie algebras.
This notion was introduced by Krichever and Novikov
\cite{KNFa} (they called it quasi-graded) and they 
constructed such an almost-grading in the higher genus and 
two point case.
To find an almost-grading in the multi-point case is more difficult.
This weaker grading is enough to  introduce and study
highest weight representations.
As a remark aside: with a special choice
of basis one has $L_1=0$ and the $L_2$ depends in a known manner on the 
genus $g$ and the number of points in $I$ and $O$ \cite{Schlce}.

The higher genus multi-point current algebra $\gb$ is the tensor
product
$\gb=\g\otimes \A$ with the Lie product
\begin{equation}
[x\otimes f,y\otimes g]=[x,y]\otimes (f\cdot g)
\end{equation}
and almost-grading 
\begin{equation}
\gb=\bigoplus_{n\in\Z}\gb_n, 
\qquad \gb_n=\g\otimes \A_n.
\end{equation}
%%%%%%%%%%%%%%%%%%%%%%%%%%%%%%%%%%%%%%%%%%%%%%%
\section{Central extensions in higher genus}
\label{S:hgcent}
The next task is to construct central extensions and
to study the question of uniqueness.
\begin{proposition}
(\cite{Schlaff})
Let $\a$ be an invariant, symmetric bilinear form
of $\g$ and $C$ a closed contour on $X$ not meeting $A$, then
\begin{equation}\label{E:cocycal}
\psi_{\alpha,C}(x\otimes f,y\otimes g) 
:=\a(x,y)\oint_C fdg
\end{equation}
is a Lie algebra two-cocycle for the current algebra $\gb$.
Hence, it defines a central extension $\gh_{\alpha, C}$.
\end{proposition}
Consequently, there  exist central extensions for $\gb$.
But contrary to the classical situation, even 
 if $\g$ is a simple Lie algebra, there will not be 
a unique nontrivial cocycle class.
If we choose essentially different contours $C$ the 
$\psi_{\a,C}$ define essentially different central extensions
 $\gh_{\alpha, C}$.

But in the classical situation we were able to extend our grading
of $\gb$ to $\gh$ by assigning a degree to the central element $t$.
This will not necessarily be true for all cocycles of the form 
\refE{cocycal}.
\begin{definition}
A 2-cocycle is called {\it local} if there exist $T_1$ and $T_2$ such that
\begin{equation}
\psi(\gb_n,\gb_m)\ne 0\implies T_2\le n+m\le T_1.
\end{equation}
\end{definition}
Given a local cocycle $\psi$ defining a central extension, then 
by setting $\deg(t)=0$ (or any other number)
the almost-grading of $\gb$ extends to the central extension.
Vice versa, if such an extension of the almost-grading exists,
the defining cocycle will be local.

We use $\Ho^2_{loc}(\gb,\C)$ to denote the subspace of 2-cocycle classes
containing a
representative which is local.
In general, the cocycles $\psi_{\alpha,C}$ are not local. But  if we choose as 
integration contour a smooth contour $C_S$ separating the points in $I$ from
the
points in $O$ and which is of winding number 1, then it can be
shown that $\psi_{\alpha,C_S}$ is a local cocycle.
In fact its values can be calculated as
\begin{equation}\label{E:csc}
\psi_{\alpha,C_S}(x\otimes f,y\otimes g) 
=\a(x,y)\oint_{C_S} fdg
=\sum_{i=1}^K\res_{P_i}(fdg).
\end{equation}
We call any such $C_S$ a separating cycle.
\begin{theorem} \cite{Schlaff}
Let $\g$ be a finite-dimensional simple Lie algebra, 
$\gb=\g\otimes A$ its (higher genus) current algebra, then
\begin{equation}
\dim \Ho^2_{loc}(\gb,\C)=1, 
\end{equation}
and a basis is given by the class of \refE{csc}, where $C_S$ is a separating
cycle and $\alpha$ is a multiple of the Cartan--Killing form.
In particular, there exists 
up to equivalence and rescaling a unique almost-graded non-trivial 
central extension of $\gb$.
\end{theorem} 
As a side-result we obtain that every local cocycle is cohomologous to a
geometric
cocycle of the form \refE{csc}.
\begin{remark}
The cocycles coming from Fock space representations and other type of
representations
are local. Hence we obtain that there exists a unique 
equivalence class (up to rescaling of the central element) of central 
extensions coming from these representations.
\end{remark}
In \cite{Schlaff} 
the more general situation of reductive Lie algebras is considered.
Recall that a finite-dimensional Lie algebra $\g$ is reductive 
if and only if it is the direct sum (as Lie algebra)
\begin{equation}
\g=\g_0\oplus\g_1\cdots\oplus\g_M, \quad \g_0 \text{ abelian },
\g_1,\g_2,\ldots,\g_M \text{ simple }.
\end{equation}

If the abelian summand $\g_0$ is missing then $\g$ is semi-simple.
In the semi-simple case it is shown that every local 2-cocycle
is cohomologous to a cocycle of the type \refE{csc} where 
$\alpha$ is an arbitrary linear combination of the individual
Cartan--Killing forms of the summands (trivially extended to the rest).
Vice versa, such cocycles are local.
In particular, we get
$\dim \Ho^2_{loc}(\gb,\C)=M$

In the reductive case we have to add another condition.
We denote by $\L$ the Lie algebra of meromorphic  vector fields on $X$ which
are
holomorphic outside of $A$.
A 2-cocycle is called $\L$-invariant if
\begin{equation}
\psi(x\otimes (e\ldot f),y\otimes g)+
\psi(x\otimes  f,y\otimes (e\ldot g))=0,\quad \forall f,g\in\A, 
\quad \forall e\in \L.
\end{equation}
Cocycles of the form \refE{cocycal} are  $\L$-invariant.
In \cite{Schlaff} it is shown that every $\L$-invariant local cocycle 
is cohomologous to \refE{csc}.
In particular
$\dim \Ho^2_{\L,loc}(\gb,\C)=M+\frac {m(m+1)}{2}$.
Here the index $\L$ denotes the classes containing cocycles as representatives
which are $\L$-invariant, and $m=\dim \g_0$. 
In the semi-simple case there is no need to pose explicitly 
$\L$-invariance as in this case in every local cocycle class there is a
unique $\L$-invariant representative.
Again, the condition of $\L$-invariance is automatic if the 
representations under consideration are representations coming from
the larger algebra of Lie algebra valued differential operators.

\begin{example}
The Lie algebra of trace-less matrices $\sln(n,\C)$ is simple. Hence, the unique
non-trivial almost-graded central extension is given 
(up to equivalence and rescaling) by the cocycle
\begin{equation}\label{E:slcoc}
\psi_1(A\otimes f,B\otimes g)=\tr(A\cdot B)\oint_{C_S}fdg.
\end{equation}
\end{example}

\begin{example}
The Lie algebra of all matrices $\gl(n,\C)$ is the direct sum
$\gl(n,C)=\mathfrak{s}(n,\C)\oplus \sln(n,C)$, where $\mathfrak{s}(n,\C)$ 
is the abelian summand of scalar
matrices. In particular $\gl(n,\C)$ is a reductive Lie algebra.
Following the general results 
\newline
$\dim \Ho^2_{\L,loc}(\gl(n,\C),\C)=2$.
A basis is given by  the elements $\psi_1$  \refE{slcoc} and 
\begin{equation}\label{E:glcoc}
\psi_2(A\otimes f,B\otimes g)=\tr(A)\cdot \tr(B)\oint_{C_S}fdg.
\end{equation}
\end{example}

%%%%%%%%%%%%%%%%%%%%%%%%%%%%%%%%%%%%%%%%%%%%%%%
\section{An example: The three-point genus zero case}\label{S:zero}
Let us consider the Riemann sphere $S^2=\P^1(\C)$ and the
set $A$ consisting of 3 points. Given any triple of 3 points there
exists always an analytic automorphism of $S^2$ mapping this triple
to $\{a,-a,\infty\}$, with $a\ne 0$.
In fact $a=1$ would suffice.
Without restriction we can take
$$
I:=\{a,-a\}, \quad  O:=\{\infty\}.
$$
Due to the symmetry of the situation it is more convenient to take a
symmetrized basis of $\A$:
\begin{equation}
A_{2k}:=(z-a)^k(z+a)^k,
\qquad
A_{2k+1}:=z(z-a)^{k}(z+a)^{k}, \qquad k\in\Z.
\end{equation}
It is shown in \cite{Schltori} that it a basis.
By more or less direct calculations one can show 
the structure equation  for the current algebra
$\gb$
\begin{equation}\label{E:three}
[x\otimes A_n, y\otimes A_m]
=\begin{cases}
[x,y]\otimes A_{n+m},& \text{$n$ or $m$ even},
\\
[x,y]\otimes A_{n+m}+a^2[x,y]\otimes A_{n+m-2},
&\text{$n$ and $m$  odd},
\end{cases} 
\end{equation}
Again $a=1$ could be set.
The reason to keep $a$ is that it can be seen that
if we vary $a$ over the affine line we obtain for $a=0$
the classical current algebra.
In particular, this family gives a deformation.
In \cite{SchlFiaaff} it was shown that this deformation
is a geometrically non-trivial deformation despite the fact 
that for $\g$ simple, $\gb$ is formally rigid, i.e. there
does not exists a non-trivial formal deformation.
This effect is peculiar to infinite dimensional Lie algebras and
is discussed in detail there.

For the central extension $\gh_{\a,S}$ of \refS{hgcent} we obtain the
defining cocycles (see \cite{SchlFiaaff}, A.13 and A.14)
\begin{equation}
\gamma(x\otimes A_n,y\otimes A_m)=
\alpha(x,y)\cdot\cins A_ndA_m,
\end{equation}
with 
\begin{equation}
\cins A_ndA_m=
\begin{cases}
\qquad -n\delta_m^{-n},&\text{$n$, $m$ even},
\\
\qquad\qquad 0,&\text{$n$, $m$ different parity},
\\
-n\delta_m^{-n}+
a^2(-n+1)\delta_m^{-n+2},&\text{$n$, $m$ odd.}
\end{cases}
\end{equation}
Of course, given a simple Lie algebra $\g$ with generators and 
structure equations the relations can  be written in these terms.
For $\sln(2,C)$ with the standard generators
$$
h:=\begin{pmatrix} 1&0\\ 0&-1 
\end{pmatrix}, \quad
e:=\begin{pmatrix} 0&1\\ 0&0 
\end{pmatrix}, \quad
f:=\begin{pmatrix} 0&0\\ 1&0 
\end{pmatrix}
$$
we set 
$e_n:=e\otimes A_n$, $n\in\Z$ and in the same way $f_n$ and $h_n$.
Recall that $\a(x,y)=\tr(x\cdot y)$.
We calculate
\begin{align}
[e_n,f_m]&=\begin{cases}  h_{n+m},& \text{$n$ or $m$ even}
  \\                    h_{n+m}+a^2h_{n+m-2},& \text{$n$ and $m$ odd}
\end{cases}
\\
[h_n,e_m]&=\begin{cases}  2e_{n+m},& \text{$n$ or $m$ even}
  \\                    2e_{n+m}+2a^2e_{n+m-2},& \text{$n$ and $m$ odd}
\end{cases}
\\
[h_n,f_m]&=\begin{cases}  -2f_{n+m},& \text{$n$ or $m$ even}
  \\                    -2f_{n+m}-2a^2f_{n+m-2},& \text{$n$ and $m$ odd}
\end{cases}
\end{align}
For the  central extension we get
\begin{equation}
[e_n,f_m]=\begin{cases}  h_{n+m}-n\delta_m^{-n},& \text{$n$ or $m$ even}
  \\                    h_{n+m}+a^2h_{n+m-2}-n\delta_m^{-n}
                         -a^2(n-1)\delta_m^{-n+2}  ,& \text{$n$ and $m$ odd}
\end{cases}
\end{equation}
The other commutators stay the same.

\newpage
\section{An example: The torus case}\label{S:torus}
%%%%%%%%%%%%%%%%%%%%%%%%%%%%%%%%%%%%%%%%%%%%%%%%%
Let $T=\C/L$ be a complex one-dimensional torus, i.e. a Riemann
surface of genus 1.
Here $L$ denotes the lattice $L=\langle 1,\tau\rangle_{\Z}$ with
$\mathrm{im}\; \tau>0$.
The field of meromorphic functions on $T$ is
generated by the doubly-periodic Weierstra\ss\ $\wp$ function and
its derivative $\wp'$ fulfilling the differential equation
\begin{equation}\label{E:diffequ}
(\wp')^2=4(\wp-e_1)(\wp-e_2)(\wp-e_3)=
4\wp^3-g_2\wp-g_3.
\end{equation}
with the $e_i$ pairwise distinct and given by  
\begin{equation}
\wp(\frac {1}{2})=e_1,\qquad
\wp(\frac {\tau}{2})=e_2,\qquad
\wp(\frac {\tau+1}{2})=e_3, \quad e_1+e_2+e_3=0.
\end{equation}
The function $\wp$ is an even meromorphic function 
with poles of order two at the
points of the lattice and holomorphic elsewhere.
The function $\wp'$ is an odd meromorphic function with poles of order 
three at the
points of the lattice  and holomorphic elsewhere.
It has zeros of order one at the points $1/2,\tau/2$ and $(1+\tau)/2$
and all its translates under the lattice.

We consider the subalgebra of functions which are 
holomorphic outside of $\bar z=\bar 0$ and $\bar z=\overline{1/2}$.
As shown in \cite{Schltori} 
a basis is given by 
\begin{equation}
A_{2k}=(\wp-e_1)^{k},
%\\
\qquad
A_{2k+1}=\frac 12 \wp'\cdot (\wp-e_1)^{k-1} \qquad k\in\Z.
%\end{aligned}
\end{equation}

See also \cite{RDS} for a similar result in the 
vector field algebra case.
The following is shown in \cite{SchlFiaaff}

\begin{equation}\label{E:structaff}
[x\otimes A_n, y\otimes A_m]
=\begin{cases}
[x,y]\otimes A_{n+m},& \text{$n$ or $m$ even},
\\
[x,y]\otimes A_{n+m}+3e_1 [x,y]\otimes A_{n+m-2}
\\
\ +(e_1-e_2)(2e_1+e_2) [x,y]\otimes A_{n+m-4},
& \text{$n$ and $m$ odd}.
\end{cases}
\end{equation}
If we let $e_1$ and $e_2$ (and hence also $e_3$) go to zero
we obtain the classical current algebra as degeneration.

The cocycle defining the central extension is given by
(\cite{SchlFiaaff}, Thm.~4.6)
\begin{equation} 
\label{E:famcent}
\gamma(x\otimes A_n,y\otimes A_m)=
\alpha(x,y)\cdot\cins A_ndA_m
\end{equation}
with 
\begin{equation}
\label{E:famcent1}
\cins A_ndA_m=
\begin{cases}
\qquad -n\delta_m^{-n},&\text{$n$, $m$ even},
\\
\qquad\qquad 0,&\text{$n$, $m$ diff. parity},
\\
-n\delta_m^{-n}+
3e_1(-n+1)\delta_m^{-n+2}+
\\
+(e_1-e_2)(2e_1+e_2)(-n+2)\delta_m^{-n+4}
,&\text{$n$, $m$ odd}.
\end{cases}
\end{equation}

\end{document}